\documentstyle[12pt]{article}

\catcode`@=11
\long
\def\@makefntext#1{\noindent #1}
\newskip\tabcentering \tabcentering=1000pt plus 1000pt minus 1000pt

\def\MCH#1#2{\setbox0=\hbox{\raise#1\hbox{#2}}\smash{\box0}}

\newcommand{\OO}{\Omega}
\newcommand{\dd}{\delta}

\newcommand{\beqs}{\begin{eqnarray*}}
\newcommand{\eeqs}{\end{eqnarray*}}
\newcommand{\beqn}{\begin{eqnarray}}
\newcommand{\eeqn} {\end{eqnarray}}
\def\pa{\partial}
\def\var{\varepsilon}

\def\beq{\begin{equation}}
\def\eeq{\end{equation}}

\begin{document}

\title{\Large THE VORTEX DYNAMICS OF A GINZBURG-LANDAU SYSTEM UNDER PINNING EFFECT
\thanks{
The first author is supported by National 973-Project and the Grant for 21th century excellent young scholars of the Education Department of China.
This work was completed when the first author was visiting
the Mathematical Department of National University of Singapore under the
research grant PR3981646. He would like to thank them for hospitality and
support. }}

\author{Huai-Yu Jian$^1$ , Xingwang Xu$^2$  \\
 \noindent{\small $^1$ Department of Mathematical Sciences, Tsinghua
University}\\
\noindent{\small Beijing 100084, P. R. China \ \ (E-mail:
hjian@math.tsinghua.edu.cn)}\\
\noindent{\small $^2$ Department of Mathematics, National University of
Singapore}\\
\noindent{\small Singapore 119260 \ \ (E-mail: matxuxw@math.nus.edu.sg)}\\}
\date{}
\maketitle

\centerline{\bf abstract}

We study the vortex dynamical behaviour of a Ginzburg-Landau (G-L) system of
related to inhomogeneous superconductors as well as to
three-dimensional superconducting thin films having variable thickness. It
is proved that the vortices are attracted by impurities or inhomogeities
in the superconducting materials. This rigorously verifies the fact
predicted recently by a few authors using a method of formal
asymptotics or approxiamate computation.
Using this fact, furthermore, we prove the strong $H^1$-convergence of
the solutions to the G-L system.


\section{Introduction and Main Results}

 Consider the solutions,
$V_\var=(V^1_\var,V^2_\var): \Omega \longrightarrow R^2$ for a
smooth bounded domain $\Omega \subset R^N$ with $N \ge 2$, to the
following initial-boundary value problem of Ginzburg-Landau
system;
 \begin{eqnarray}  \label{eq1.1}
\left\{
\begin{array}{clll}
\displaystyle{\frac{\partial V_{\varepsilon}}{\partial t}} & = & \Delta V_{\varepsilon} +\nabla \omega
\nabla V_{\varepsilon}+AV_{\varepsilon}+ \displaystyle{\frac{BV_{\varepsilon}}{\varepsilon^2}}(1-|V_{\varepsilon}|^2)\,\,\,\, & %
\mbox{in $\OO \times (0, \infty )$} \\
V_\var(x,t) & = & g(x) & \mbox{on $\pa\OO \times (0, \infty )$} \\
V_\var(x,0) & = & V^0_{\varepsilon}(x) & \mbox {in $\OO$},
\end{array}
\right.
\end{eqnarray}
where $\omega ,A$ and $B$ are given smooth functions defined on ${\bar
\Omega }.$ Equation (\ref{eq1.1}) is a simple model which simulates
inhomogeneous type II superconducting materials as well as three-dimensional
superconducting thin films having variable thickness. In the inhomogeneous
materials, the equilibrium density of superconducting electrons is not a
constant, but a positive and smooth function on $\bar \Omega $. Letting $%
a(x) $ denote the density function and neglecting the magnectic field, one
obtains that
\begin{equation}  \label{eq1.2}
\displaystyle{\frac{\partial u}{\partial t}}- \Delta u = \displaystyle{\frac{%
u}{\varepsilon^2}}(a(x)-|u|^2).
\end{equation}
By setting $V(x,t)= u(x,t)/ \sqrt{a(x)}$, equation (\ref{eq1.2}) reduces to
the equation in (\ref{eq1.1}) with
\[
\omega (x)=\ln a(x), \ \ A(x)=\frac{\Delta \sqrt{a(x)}}{\sqrt{a(x)}} \ \ and
\ \ B(x)=a(x).
\]
We refer to \cite{cr}, \cite{cdg}, \cite{r} and the references therein for
the detailed discussion of the motivation and physical background for
equation (\ref{eq1.2}) and its more general form involving magnectic field
and electric field. There are several theoretical results on the static case
of (\ref{eq1.2}); see, for example, \cite{r}, \cite{as}, \cite{dl}, \cite{j1}
for the minimizers of the functional associated with the static case and
\cite{jw1} for the state-solutions.

Another model which also reduces to the equation in (\ref{eq1.1}) is the
following:
\begin{equation}  \label{eq1.3}
\displaystyle{\frac{\partial u}{\partial t}}=\displaystyle{\frac{1}{a(x)}}%
div (a(x)\nabla u)+ \displaystyle{\frac{u}{\varepsilon^2}}(1-|u|^2).
\end{equation}
This equation is related to three-dimensional superconducting thin films
having variable thickness. Let $\Omega \times (-\delta a(x) , \delta a(x))$
be the domain occupied by these materials, where $\Omega\subset R^2.$ Then
this superconducting film was modeled as two-dimensional objects by equation
(\ref{eq1.3}) in \cite{dg}, \cite{cdg1}. We refer the reader to \cite{as},
\cite{dg} and \cite{cdg1} for the study of the minimizers of a functional
associated with the static case of (\ref{eq1.3}). Obviously, (\ref{eq1.3})
is the special form of the equation in (\ref{eq1.1}) with
\[
\omega (x)=\ln a(x), \ \ A\equiv 0 \ \ and \ \ B\equiv 1.
\]
The connection between the steady solutions for (\ref{eq1.1}) and the
self-similar solutions for harmonic maps was studied in \cite{jw2} recently.

Physically, the points at which a solution to problem (\ref{eq1.1}) equals
to zero are called vortices. In the case of $N=2, \omega (x)\equiv
A(x)\equiv 0$ and $B(x)\equiv 1,$ the vortex dynamics was studied previously
for the steady equations by Bethuel, Brezis and H\'{e}lein \cite{bbh}. (For the minimum solution, see \cite{s}).
Furthermore, Lin \cite{l2}, independently, Jerrard and Soner \cite{js1}
and \cite{js},
studied the dynamical law for the vortices of $u_{\varepsilon }(x,
t)=V_{\varepsilon } (x, |\ln \varepsilon |t),$ where $V_{\varepsilon } (x,
t) $ solves the initial-boundary value problem (\ref{eq1.1})
under the same case. Their
dynamical law is described by an ODE, ${\frac {d}{dt}}y(t)=-\nabla W(y(t))$
where $W$ is some known function related to the domain and the boundary
condition
and called as the
 renormalized energy functional associated with the steady problem
\cite{bbh}, \cite{l2}or \cite{js}. The results in \cite{l2} and in \cite{js}
were generalized to the Neumann boundary condition
by Lin in \cite{l3}.

However, there are few results for the vortex dynamics in the orginal time
(not scaling by the time factor $|\ln \varepsilon |$), especially for
equations (\ref{eq1.2}) and (\ref{eq1.3}), not to speak of for (1.1).
Up to our limit knowledge, one
can only locate one result for equation (1.3) by Lin in \cite{l3}. He proved that as $%
\varepsilon \to 0 ,$ under some suitable assumptions
 on the initial and boundary data, the
solutions, $u_\var (x,t) $, of the Dirichlet initial-boundary valued problem
for equation (\ref{eq1.3}), subconverge
 in $H^1_{loc}(\bar\Omega \setminus \{y_1(t), y_2(t), \cdots, y_k(t)\}),$
where
 the
functions $y_j(t):[0, T )\to \Omega\subset R^2$ satisfy the following ODE:
\begin{eqnarray}  \label{eq1.4}
\left\{
\begin{array}{cll}
{\frac{d }{d t}}y_j(t) & = & - a^{-1}(y_j(t)))\nabla a(y_j(t)) , \ \ 0\leq
t\leq T , \\
y_j (0) & = & b_j
\end{array}
\right.
\end{eqnarray}
Here $k, d_j , d$ are some constants related to the initial data, while $T$
is chosen so that $a_j(t)$ will stay inside $\Omega$ and $y_l(t)\not%
{=}y_j(t) $ for all $0\leq t\leq T$ and for all $j\not{= }l,$ $j,l=1,2,
\cdots, k$. See Theorem 1.1 in \cite{l3} for the details. The first author
in \cite{j2} proved that $T,$ in fact, is $+\infty $ and each solution $%
y_j(t)$ of (\ref{eq1.4}) converges to a critical point of $a(x)$ as $t\to
\infty$ as long as $a(x)$ is analytic at its critical points and $b_j \in
\Omega_j \subset \subset \Omega $ satisfies
$
\min_{x\in \partial \Omega _j }a(x)> a(b_j)
$
for some domain $\Omega _j .$

 In recent paper [18], the authors started studying the vortex dynamics of
 equation (1.2) with $N=2$. They proved that the vortices are attracted by the local minimum points
 of $a(x)$
 and the vortex dynamics is described by equation (1.4). Moreover, the authors conjectured that similar results
  should be true for equation (1.3).

In this paper, we will verify this conjecture. In fact, we will prove that
all results for equation (1.2) in [18] are also true for problem (1.1)
 (see
Theorem 1 and the first part of Theorem 2 below). In particular, we will prove that for
most  sufficient large $t$, under some suitable conditions, all the
vortices of problem (\ref{eq1.1}) are pinned together to the local minimum points
of $\omega (x)$ in $\Omega $ as $\varepsilon \to 0 .$ This result was
observed by Chapman and Richardson in \cite{cr} for equation (\ref{eq1.2})
and Du and Gunzburger in \cite{dg} for equation (\ref{eq1.3}). They used a
matched asymptotic method or approximate computation method to predict that
vortices for equation (\ref{eq1.2}) or (\ref{eq1.3}) (in fact, for a more
complicated equation involving magnetic field and electric field), are
attracted to the the local minimum points of $a(x).$ Our results in this paper
will show that their predictions are correct. See Remark 1  below.

As our second goal, we will study the strong $H^1$-convergence of solutions
to problem (1.1). Although our strong convergence result, the second part of
Theorem 2 below, can be viewed as a generalization of Theorem 1.1 in [16],
our method to prove it is based on the vortex convergence,
the first part of Theorem 2, and is completely different from the arguments in [16].

\vskip 0.3cm

To state our main results, we need the following assumptions:

${\bf (A_1)}$ \ \ {\it $g:\partial\Omega\longrightarrow R^2$ is smooth, $%
|g(x)| \equiv 1 $ on $\partial\Omega ;$  }

${\bf (A_2)}$ \ \ {\it $\omega \in C^{2, \alpha }(\bar \Omega ), A, B\in
C^{\alpha} (\bar \Omega ) \ \ (\alpha >0),$ $\omega(x)>0$ and $B(x)>0$ }

{\it \ \ \ \ \ \ \ \ \ \ for all $x\in \bar \Omega;$}

${\bf (A_3)}$ \ \ {\it the initial data $V^0_{\varepsilon}\in C^2(\bar
\Omega;R^2)$ \ \ $(\varepsilon >0)$ satisfies $V^0_{\varepsilon }(x)= g(x) $
}

{\it \ \ \ \ \ \ \ \ \ \ on $\partial\Omega $ and
\[
||V^0_{\varepsilon }||_{C(\bar \Omega )}\leq K,\ \ \int _{\Omega} \rho ^2(x)
[|\nabla V^0_{\varepsilon}|^2 +{\frac {1}{2\varepsilon ^2}}%
(|V^0_{\varepsilon}|^2-1)^2]dx \leq K
\]
\noindent for a constant $K$ (independent of $\varepsilon $ ) and some $m$
distinct points $b_1 ,b_2 , \cdots , b_m$ in $\Omega , $ where $\rho(x)= min
\{|x-b_j|, j=1,2, \cdots , m\}.$ }

${\bf (A_4)}$ \ \ {\it For each $j$, there exists some Lipschitz domain $G_j$
such that $\ \ b_j \in G_j \subset \subset \Omega, \min_{x\in \partial G_j
}\omega(x)> \omega(b_j), j=1, \cdots , m$.}

\vskip 0.3cm As a start point, consider the ordinary differential system
\begin{eqnarray}  \label{eq1.5}
\left\{
\begin{array}{cll}
{\frac{d }{d t}}y_j(t) & = & - \nabla \omega (y_j(t)) , \ \ 0\leq t < \infty ,
\\
y_j (0) & = &b_j
\end{array}
\right.
\end{eqnarray}
\noindent where $j=1,2, \cdots, m,$  and
$\nabla \omega $ is the gradient of the function $\omega $ with respect to $%
x=(x_1,\cdots , x_N)\in R^N.$ It is this system that will play an important
rule in the course of the proof of our main results. As preliminary, we will
generalize the resuls for (1.4) in [17] to  system (1.5)
in next section under an extra condition

${\bf (A_5) }$ \ \ {\it either $\omega$ has only nondegenerate critical points
in $\Omega $ or it is analytic in some neighborhood of its critical points.}

Particularly, the existence and uniqueness of global solutions to (1.5)
will be guaranteed by conditions $(A_2)$ and $(A_4) $(see Lemmas 4 and 5 below).

\vskip 0.3cm

As the main results, we will first prove the following compactness for the
solutions to problem (\ref{eq1.1}) for all $N\geq 2$ which is the
generalization of Theorem 1.1 in \cite{l3} and Theorem 1.2 in [18]
for $N=2$.

{\bf Theorem 1}
 {\it Suppose that the hypotheses $(A_1),(A_2),(A_3)$ and $(A_4)$
are satisfied. Let $y_j(t)$ be solutions to problem (\ref{eq1.5})
$(1\leq j\leq m)$ and set
\[
\Omega (\omega )=\bar{\Omega}\times (0,\infty )\backslash \cup
_{j=1}^m\{(x,t):x=y_j(t),0<t<\infty \}.
\]
Then there is a positive constant $\varepsilon _0$ (depending only on the
infimum of $B$ and the superemum of $|A|$) such that the set $%
\{V_\varepsilon :\ \ \varepsilon \in (0,\varepsilon _0)\}$ of the classical
solutions to problem (\ref{eq1.1}) is bounded in $H_{loc}^1(\Omega (\omega
)).$ Moreover, given any sequence $\varepsilon _n\downarrow 0,$ there exists
a subsequence (still denoted by itself) such that
$V_{\varepsilon
_n}\longrightarrow V$
weakly in $H_{loc}^1(\Omega (\omega ))$ with $V$ satisfying $|V(x,t)|=1$
a.e. in $\Omega (\omega )$, $V=g$ on $\partial \Omega \times (0,\infty ),$
and the equation
\[
\displaystyle {\frac{\partial V}{\partial t}=\Delta V+V|\nabla V|^2+\nabla
\omega \nabla V\ \ in\ \ D^{\prime }(\Omega (\omega ))}.
\] }

\vskip 0.3cm We will use this compactness result, covering arguments and elliptic
estimates
to prove the following vortex convergence.

{\bf Theorem 2}
 {\it Suppose that $N=2$ and all the
hypotheses in Theorem 1 under this case
are fulfilled.
Let $V_{\varepsilon}$ be a classical solution to problem (1.1) for each
$\varepsilon >0$. Then
for
 any $\delta \in (0,\frac 14)$ and any interval $I\subset (0, \infty )$
 with $|I|>0$,
 one can find  $t\in I$ and  $\varepsilon_1 > 0$ such that the following two
 conclusions  hold true for all $\varepsilon \in (0, \varepsilon _1)$:

 {\bf (i)}
\[
\{x\in \bar{\Omega}:|V_{\varepsilon }(x,t)|\leq \frac 12\}\subset \cup
_{j=1}^mB_\delta (y_j(t))
\]
and
\[
\displaystyle { ||V_{\varepsilon }}(x,t)|-1|\leq C(\delta ,t)\varepsilon
^{\frac 12},\ \ \forall x\in \bar{\Omega}\setminus \cup _{j=1}^mB_\delta
(y_j(t)) := \Omega (\omega^\delta _t);
\]

{\bf (ii)}  \ \  if , furthermore, $A(x)\leq B(x) $ for all $x\in \bar{ \Omega}$,
then
\begin{equation}  \label{eq1.6}
 |||V_{\varepsilon }(\cdot,t)|-1||_{H^1(
\Omega (\omega^\delta _t))} +
\varepsilon^{-1}|||V_{\varepsilon }(\cdot,t)|-1||_{L^2(
\Omega (\omega^\delta _t))}
\leq C(\delta ,t)\varepsilon^{\frac 14}
\end{equation}
and the convergence,
$V(\cdot , t)_{\varepsilon
_n}\longrightarrow V(\cdot , t) , $ in Theorem 1, is strong convergence in
$H^1(\Omega (\omega^\delta _t))$ as $n\to \infty$.}

{\bf Corollary 3}
 {\it Besides the assumptions in Theorem \ref{th2}, we further
assume that for some $t$, the set

$\{\partial _tV_\varepsilon (x,t),\varepsilon \in (0,\varepsilon _0)\}$ is a
bounded set in $L_{loc}^2({\bar{\Omega}}\backslash \{y_1(t),\cdots
,y_m(t)\}).$
Then for any $\delta \in (0,\frac 14)$, there is a constant $\varepsilon
_1>0 $ such that
conclusions (i) and (ii) of Theorem 2 hold true
for all $\varepsilon \in (0,\varepsilon _1)$. }

\vskip 0.3cm

The organization of this paper is as follows. As we have mentioned before,
we will follow the arguments in [17] to  study the system (\ref{eq1.5}) and generalize the
main result in \cite{j2} for our use later in next section. Section three
will be devoted to the proof of Theorem 1. Finally in section
four, we will finish the proof of Theorem 2 and its corollary 3.

Before we are going to the detail proof of our results, we would like to
make the following remarks.

{\bf Remark 1}
 {\it It is worth to point out that Theorem 2 (Corollary
3),together with Lemmas 4 and 5 below,  imply that for most large $t$, all the vortices of $V_\varepsilon
(x,t)$ for (\ref{eq1.1}), as $\varepsilon \to 0$, move toward and
eventurally pin together at the critical points of $\omega (x)$ in $\Omega $
if $(A_1)-(A_5)$ are satisfied. In
particular, if
 each $\bar{G_j}$ contains no other
critical points of $\omega (x)$ than its local minimum points, then all the
vortices of $V_\varepsilon (x,t)$ are pinned to the minimum points.}

{\bf Remark 2} {\it Obviously, the unique solution of (\ref{eq1.5}) is $%
y_j(t)\equiv p_j\ \ (j=1,2,\cdots ,m)$ for all $t\in [0,+\infty )$ if $%
\nabla \omega (p_j)=0.$ This shows that if $\omega \equiv constant$, the
vortices for (\ref{eq1.1}) do not move in any finite time interval. In fact,
in this case, formal analyses indicate that, if initial $V_0^\varepsilon $
has isolated vortices, then these vortices move with velocities of the order
of $|\ln \varepsilon |^{-1}$. This prediction was proved rigorously in \cite
{l2} and \cite{js}.}

\vskip 0.3cm

{\bf Convention:} \ \ {\it Throughout this paper, we use the letter $C$ to
denote various constants independent of $\varepsilon$ but maybe depending on
$\Omega, \omega (x) , A(x),B(x), g, K$ and other known constants.}

\section{ Preliminaries}

\setcounter{equation}{00}

Observing that ordinary differential system (1.5) is nothing but (1.4)
with $ln a(x)$ replaced by $\omega (x)$, we can apply the result and method
for (1.4) in [17] to (1.5).

First, following the arguments from (2.1) to (2.3) in [17], we have

{\bf Lemma 4}
 {\it Suppose $\omega $ satisfies $(A_2)$ and $(A_4)$.
Then (1.5) has a unique
$C^3$ solution $(y_1,y_2,\cdots ,y_m):[0,\infty)\rightarrow {R}^{mN}$.
Furthermore,
if $b_j\not= b_l$ for all $j\not=l$, then

\begin{enumerate}

\item[a.]  for all $t\in [0,\infty )$ and $1\le j\not=l\le m$, $y_j(t)\not
=y_l(t)$;

\item[b.]  for all $t\in [0,\infty )$ and $1\le j\le m$, $y_j(t)\in G_j$.
\end{enumerate}
}
 {\bf Lemma 5}
{\it  Assume that $\omega $ satisfies $(A_2)$ and $(A_4)$ as well as
$(A_5)$. Then for each $j=1,2,\cdots ,m$, there exists $b_j\in {\bar{G}_j}$
which are critical point of $\omega $ such that $y_j(t)\rightarrow b_j$ as $%
t\to +\infty $.  }

\noindent{\bf Proof.}  Repeating the proof of Theorem 1.2 in [17]
with $ln a(x)$ replaced by $\omega (x)$, one easily obtains the desired results.
Here we would like to point out two facts in order for the reader easy to
follow the arguments.

1. To prove $y(t)\to b$ as $t\to \infty$ under the conditions that
   $\nabla \omega(b) =0 $ and
$\omega (y(t))\downarrow \omega (b),$
  we may assume that $\omega (y(t))>\omega (b)$ for
all $t\in (0,\infty ).$ Otherwise, it is easy to find a $t_0>0$ such that $%
\omega (y(t))=\omega (b)$ for all $t\in (t_0,\infty ).$ This yields
\[
0=\nabla \omega (y(t))\cdot y^{\prime }(t)=-|y^{\prime }(t)|^2
\]
by (1.5). Hence $y(t)=b$ for all $t\in (t_0,\infty ).$

2. The hypothesis that $\omega $ has only nondegenerate critical points $b$
implies that
\begin{equation} \label{eq2.1}
|\nabla \omega(x)|\geq \theta_1 | \omega(x) - \omega(b)|^{\frac{1}{2}}
\ \
\forall x\in B_{\theta_2}(b)
\end{equation}
for some positive constants $\theta_1$ and $\theta_2 .$

In fact, Since           $\nabla \omega(b) =0 $   and
$det (\nabla ^2 \omega (b))\not{=}0,$  we have
\begin{eqnarray*}
|\nabla \omega(x) | &=& |\nabla^2\omega(b)\cdot (x-b) | +o(|x-b|)| \\
&=& \{(x-b)^{\top}\cdot (\nabla^2 \omega (b))^{\top} \nabla^2 \omega
(b)\cdot (x-b)\}^{\frac{1}{2}}+o(|x-b|) \\
&\geq & \lambda_0 |x-b| +o(|x-b|)
\end{eqnarray*}
and
\begin{eqnarray*}
| \omega(x)-\omega (b)| &=& |(x-b)^{\top}\cdot \nabla^2 \omega (b)\cdot
(x-b) | +o(|x-b|^2) \\
&\leq & max \{|\lambda_1|, |\lambda_2|, \cdots , |\lambda_n|\} |x-b|^2
+o(|x-b|^2),
\end{eqnarray*}
where $\lambda_0$ is the minimum eigenvalue of $(\nabla^2 \omega (b))^{\top}
\nabla^2\omega (b)$ and $\lambda_1, \lambda_2, \cdots , \lambda_n$ are the
eigenvalues of $\nabla ^2\omega (b) .$ Thus  the hypothesis
$det (\nabla ^2 \omega
(b))\not{=}0$ guarantees that (2.1) is true.
 \vskip 0.3cm

{\bf Corollary 6}
 {\it In Lemma 5, if $\Omega $ is replaced by $R^N$ and $%
(A_4)$ replaced by the following assumption:
\begin{equation}
|x|\leq \Psi (\omega (x))  \label{eq2.2}
\end{equation}
and
\begin{equation}
|\nabla \omega (x)|\leq \phi (|x|),  \label{eq2.3}
\end{equation}
where $\Psi :R\to R^{+}$ is a increasing function and $\phi $ is a positive
function satisfying
\begin{equation}
\int_\alpha ^{+\infty }\phi ^{-1}(t)dt=+\infty  \label{eq2.4}
\end{equation}
for some $\alpha >0$, then all the conclusion of Lemma 5 is also
true.            }

{\bf Proof. } Using (2.3) and (2.4)and
applying Wintner's theorem (Theorem 2.5 of Chapter 1 in \cite{zd}) to system
(1.5), one can conclude that (1.5) has a unique solution in $[0,+\infty )$.
Moreover, by (2.2), one can see that the solution is
bounded in $R^N$ uniformly for $t\in [0,\infty )$. The
remaining part of the proof is the same as in [17].

\section{Proof of Theorem 1}

 \setcounter{equation}{00}

 Throughtout this section, we assume $(A_1), (A_2), (A_3)$ and $%
(A_4),$ although some conclusions below need only parts of these
assumptions.

{\bf Lemma 8.}
 {\it Let $V_{\varepsilon}$ be classical solutions to
problem (\ref{eq1.1}). Then there exists an $\varepsilon _0>0$ such that for
all $\varepsilon \in (0,\varepsilon _0)$
\begin{equation}
|V_{\varepsilon}(x,t)|^2\leq C,\ \ \forall (x,t)\in \bar{\Omega}\times
[0,\infty )  \label{eq3.1}
\end{equation}
and
\begin{equation}
|\nabla V_{\varepsilon}(x,t)|^2+|\frac{\partial V_{\varepsilon}}{\partial t}|
\leq \frac C{\varepsilon ^2},\ \ \forall (x,t)\in \bar{\Omega}\times
[\varepsilon ^2,\infty ).  \label{eq3.2}
\end{equation}}

{\bf Proof.}  Denote
\begin{equation}  \label{eq3.3}
\varepsilon _0=\displaystyle{[\frac{\inf_{x\in \Omega }B(x)}{1+\sup_{x\in
\Omega } |A(x)|}]^{\frac{1}{2}}}.
\end{equation}
Let $W =|V_\var|^2$. Dropping the subscript $\varepsilon$, we see that the
equaton in (\ref{eq1.1}) reads as
\begin{equation}  \label{eq3.4}
\partial_t W =\Delta W +\nabla \omega \nabla W -2|\nabla V|^2 +\frac{2BW}{%
\varepsilon^2} (1-W)+2AW .
\end{equation}

If (3.1) were not true for some $\varepsilon \in (0,
\varepsilon _0)$, we could use $(A_1)$ and $(A_3)$, and employ the usual
arguements for maximum principle to find a point $(x_{\varepsilon
},t_{\varepsilon }) \in \Omega \times (0, \infty )$ (for each $\varepsilon $
) at which
\begin{equation}  \label{eq3.5}
W>2, \ \ \nabla W=0 , \ \ \partial _t W \geq 0 \ \ and \ \ \Delta W\leq 0.
\end{equation}
Moreover, (\ref{eq3.4}) gives us
$\partial _t W \leq 2(A-\frac
{B}{\varepsilon ^2})W$ at $(x_{\varepsilon },t_{\varepsilon }).$
This yields
a contradiction as $\varepsilon \in (0, \varepsilon _0).$

By a scaling arguement, considering the equation for $U_\var
(x,t)=V_\var (\varepsilon x , \varepsilon ^2 t)$ and using (\ref{eq1.1}) and
standard local parabolic estimates, we immediately obtain (\ref{eq3.2}).

\noindent{\bf Corollary 8} \ \ {\it
 Let $V_{\varepsilon}$ be classical solutions to
problem (1.1). If $A(x)\leq B(x)$ for all $x\in \bar{\Omega} ,$
then $
|V_{\varepsilon}(x,t)|^2\leq 1+\varepsilon^2$
for all $(x,t)\in \bar{\Omega}\times
[0,\infty )$ and for all $\varepsilon >0$.  }

{\bf Proof.} It is obvious from the proof of Lemma 7,
observing that $W>1+\var ^2$ and so
$\partial _t W \leq 2(A-B)W$ at $(x_{\varepsilon },t_{\varepsilon }).$

\vskip 0.3cm

 Set
\begin{equation}  \label{eq3.6}
E_\var (V(x,t))=\frac {e^{\omega (x)}}{2} [ |\nabla V(x,t)|^2 +\frac{B(x)}{%
2\varepsilon^2}(1-|V(x,t)|^2)^2].
\end{equation}

\vskip 0.3cm

\noindent{\bf Lemma 9}\ \
 {\it Let $V_{\varepsilon}$ be solutions to (\ref{eq1.1}).
Then for any $T>0,$ there exist two positive constants $C(T)$ and $\sigma
(T) $ (both depending on T) such that for all $\varepsilon \in
(0,\varepsilon _0),$ all $\delta \in (0,\sigma (T))$ and all $t\in [0,T],$
one has
\[
B_\delta (y_l(t))\subset \Omega ,\ \ B_\delta (y_l(t))\cap B_\delta
(y_j(t))=\emptyset \ \ for\ \ l\not{=}j
\]
and
\[
\int_0^T\int_{\Omega \backslash \cup _{j=1}^mB_\delta (y_j(t))}|\frac{%
\partial V_{\varepsilon}}{\partial t}|^2dxdt+\sup_{0\leq t\leq
T}\int_{\Omega \backslash \cup _{j=1}^mB_\delta
(y_j(t))}E_{\varepsilon}(V_{\varepsilon})dx\leq \delta ^{-2}C(T).
\]               }

{\bf Proof. }For each $T>0,$ by Lemma 4 we can find a $
\sigma =\sigma (T)>0$ such that
\begin{equation}
\sigma \leq \frac 2{1+\sup_{x\in \Omega }|\nabla \ln B|}  \label{eq3.7}
\end{equation}
and for all $t\in [0,T],$
\begin{equation}
\min_{1\leq l,j\leq m}\{dist(y_j(t),\partial \Omega ),\ \ |y_j(t)-y_l(t)|\ \
for\ \ l\not{=}j\}\geq 4\sigma .  \label{eq3.8}
\end{equation}
Motivated by a method in \cite{js}, we choose a
smooth monotone function $\phi :[0,\infty )\longrightarrow [0,\infty )$ such
that
\begin{equation}
\phi (r)=\left\{
\begin{array}{cc}
r^2, & if\ \ r\leq \sigma \\
\sigma ^2, & if\ \ r\geq 2\sigma .
\end{array}
\right.  \label{eq3.9}
\end{equation}

 Let
\[
\rho (x,t)=\min_{1\leq j\leq m}|x-y_j(t)|.
\]
It follows easily from (\ref{eq3.8}) that $\phi (\rho (x, t))$ is smooth in $%
x$ as well as in $t$ for all $(x, t)\in \bar \Omega \times [0, T].$ Dropping
the subscript $\varepsilon ,$ applying integration by parts, noting $%
\partial _t V=\partial _t g =0$ on $\partial \Omega \times (0, \infty )$ and
using (\ref{eq3.6}) and the equation in (\ref{eq1.1}), we obtain
\begin{eqnarray}  \label{eq3.10}
& \ \ & \frac{d }{d t}\int_{\Omega } \phi (\rho (x, t)) E(V)dx  \nonumber \\
& = & \int_{\Omega }\frac {d\phi (\rho )}{d t} E( V) +\int_{\Omega } \phi
(\rho )e^{\omega }[\nabla V\nabla V_t- \varepsilon ^{-2}B V(1-|V|^2)V_t]
\nonumber \\
& = & \int_{\Omega } \frac {d\phi (\rho )}{d t}E(V)-\int _{\Omega }V_t\nabla
V \nabla \phi (\rho )e^{\omega } -\int _{\Omega }\phi (\rho )[\nabla
(e^{\omega }\nabla V)  \nonumber \\
& \ \ & +\varepsilon ^{-2}Be^{\omega } V(1-|V|^2)]V_t  \nonumber \\
& = & \int_{\Omega } \frac {d\phi (\rho )}{d t}E(V)-\int _{\Omega }e^{\omega
}|V_t|^2 \phi (\rho )+ \int_{\Omega }e^{\omega } V_t A V \phi (\rho ) -\int
_{\Omega }e^{\omega } V_t\nabla V\nabla \phi (\rho )  \nonumber \\
& \leq & \int_{\Omega } \frac {d\phi (\rho )}{d t}E(V)- \frac {1}{2}\int
_{\Omega }\phi (\rho )e^{\omega }|V_t|^2+ \frac{1}{2} \int _{\Omega }\phi
(\rho ) e^{\omega }|AV|^2 -\int _{\Omega }e^{\omega } V_t\nabla V\nabla \phi
(\rho )  \nonumber \\
& \leq & C -\frac {1}{2}\int _{\Omega }\phi (\rho )e^{\omega }|V_t|^2+
\int_{\Omega }[ \frac {d\phi (\rho ))}{d t}E(V)- e^{\omega }V_t\nabla
V\nabla \phi (\rho )]dx ,
\end{eqnarray}
where we have used (\ref{eq3.1}) and (\ref{eq3.9}) and $C$ depends only on $%
K, \sigma , ||A||_{C(\Omega )}$ and $||\omega ||_{C(\Omega )}$.

 One can use the notation $U_i = \frac {\partial U }{\partial x_i}$
and the summation convention to compute
\begin{eqnarray*}
e^{\omega} V_t\nabla V\nabla \phi (\rho ) &=& e^{\omega }\nabla V \nabla
\phi [ \Delta V+ \nabla \omega \nabla V+A V + \varepsilon ^{-2}B V (1-|V|^2)]
\\
&=& e^{\omega }\phi _i \{ - (4\varepsilon ^2)^{-1}[B(1-|V|^2)^2]_i +
(4\varepsilon ^2)^{-1}(1-|V|^2)^2 B_i \\
& \ \ & -2^{-1}(|\nabla V|^2)_i +[V_i V_j]_j+ V_i\omega _j V_j +AVV_i\}.
\end{eqnarray*}

 By virtue of this equality, integration by parts and the fact of $%
\nabla ^k \phi (\rho )=0$ on $\partial \Omega$ (see (\ref{eq3.8}) and (\ref
{eq3.9})), one gets that
\begin{eqnarray}  \label{eq3.11}
&\ \ & \int _{\Omega }e^{\omega } V_t\nabla V\nabla \phi (\rho )dx  \nonumber
\\
&=& \int_{\Omega }[\Delta \phi E(V)+ (4\varepsilon ^2)^{-1}e^{\omega
}(1-|V|^2)^2\nabla B \nabla \phi  \nonumber \\
& \ \ &+\nabla \phi \nabla \omega E(V)-V_iV_j(e^{\omega }\phi _i)_j
+e^{\omega }(V_iV_j\omega _j+AVV_i)\phi _i]dx.
\end{eqnarray}

Combing (\ref{eq3.10}) and (\ref{eq3.11}) yields

\begin{eqnarray}  \label{eq3.12}
\frac{d }{d t}\int_{\Omega } \phi (\rho (x, t)) E(V)dx &\leq & C -\frac
{1}{2}\int _{\Omega }\phi (\rho )e^{\omega }|V_t|^2  \nonumber \\
& +& \int _{\Omega } \{E(V) [(\phi (\rho )_t-\nabla \omega \nabla \phi ]
\nonumber \\
& +& I_1(V)\}dx +I_2(V) ,
\end{eqnarray}
where
\begin{equation}  \label{eq3.13}
I_1(V)=e^{\omega } V_iV_j\phi _{ij}-\Delta \phi E(V)- (4\varepsilon
^2)^{-1}e^{\omega }(1-|V|^2)^2\nabla B \nabla \phi
\end{equation}
and
\begin{eqnarray}  \label{eq3.14}
I_2(V)&=& -\int _{\Omega } e^{\omega }A V \nabla V \nabla \phi  \nonumber \\
&\leq & \int _{\Omega }[|AV|^2\frac {|\nabla \phi |^2}{\phi }+ \phi |\nabla
V|^2 ]e^{\omega }  \nonumber \\
&\leq & C [1+\int_\OO \phi |\nabla V|^2 ].
\end{eqnarray}
Here we have used (\ref{eq3.1}) and (\ref{eq3.9}).

 If $\rho (x,t)\geq \sigma ,$ by (\ref{eq3.9}) one has
\begin{equation}  \label{eq3.15}
E(V)|\phi (\rho )_t-\nabla \omega \nabla \phi | + |I_1(V)| \leq C \phi (\rho
)E(V).
\end{equation}
If $\rho (x,t)\leq \sigma ,$ on the other hand, then $\phi (\rho (x,t))=
|x-y_l(t)|^2$ for some $l$. Hence $\phi _{ij}=\delta _{ij} $ and
\begin{equation}  \label{eq3.16}
I_1(V)= e^{\omega } (1-\frac{N}{2})|\nabla V|^2 - \frac{e^{\omega
}B(1-|V|^2)^2}{4 \varepsilon ^2}[N+ 2(x-y_l(t))\cdot\nabla \ln B ] \leq 0
\end{equation}
by (\ref{eq3.7}). Moreover, using (\ref{eq1.5}) we have
\begin{eqnarray}  \label{eq3.17}
\phi (\rho )_t-\nabla \omega \nabla \phi &=& 2(x-y_l(t))(\nabla \omega
(y_l(t))- \nabla \omega (x))  \nonumber \\
&\leq & 2|x-y_l(t)|^2||\omega ||_{C^2(\bar \Omega )}  \nonumber \\
&=& 2||\omega ||_{C^2(\bar \Omega )} \phi (\rho (x,t)).
\end{eqnarray}
Combing (\ref{eq3.12})-(\ref{eq3.17}), we obtain
\begin{eqnarray*}
& \ \ & \frac{d }{d t}\int_{\Omega } \phi (\rho (x, t)) E_\var (V_\var)dx +
\frac {1}{2}\int _{\Omega }\phi (\rho (x,t) )e^{\omega }|V_t|^2 \\
&\leq & C [1+\int_{\Omega } \phi (\rho (x, t)) E_\var (V_\var )dx]
\end{eqnarray*}
for all $t\in [0, T].$ Hence, by Gronwall's inequality and $(A_3),$ we
deduce that
\[
\int_{\Omega } \phi (\rho (x, t)) E_\var (V_\var)dx + \frac {1}{2}\int
_0^t
\int _{\Omega }E^{c(t-s)}\phi (\rho (x,t) )e^{\omega } |\partial _t V_\var
|^2dxds\leq C .
\]
This result, together with the fact that $\phi (\rho (x,t))\geq \delta ^2$
for all $t\in [0, T],$ all $x\in (\bar \Omega \backslash \cup_{l=1}^m B_\dd
(y_l(t))$ and any $\delta \in (0, \sigma (T)),$ immediately implies the
conclusion of Lemma 9.

\vskip 0.3cm

{\bf Proof of Theorem 1}: Recall the $G_j$ in $(A_4)$ and let
\begin{equation}  \label{eq3.18}
\delta _0 =\min_{1\le j\leq m} dist (G_j , \partial \Omega ).
\end{equation}
Then by Lemma 4, we have
\begin{equation}  \label{3.19}
\min_{1\leq j\leq m} dist (y_j(t), \partial \Omega )\geq \delta _0, \forall
t\in [0, \infty ).
\end{equation}
Set
\begin{equation}  \label{eq3.20}
\Omega _r(t)= \Omega \backslash \cup_{j=1}^m B_r(y_j(t)).
\end{equation}
For any $T>0$ and any $\delta \in (0, \delta _0),$ it follows from Lemma 9 and $(A_2)$ that
\begin{equation}  \label{eq3.21}
\int_0^T \int _{\Omega _{\delta /4}(t)}|\frac {\partial V_\var }{\partial
t}|^2dxdt +\sup_{0\leq t\leq T} \int _{\Omega _{\delta /4 }(t)}[|\nabla
V_\var |^2 +\frac{1}{\varepsilon ^2}(1-|V_\var |^2)^2]dx \leq C(\delta ,T)
\end{equation}
for all $\varepsilon > 0.$ This shows that the set $\{V_\var : \ \
\varepsilon >0 \}$ is bounded in $H_{loc}^1(\Omega (\omega )).$

Using (\ref{eq3.21}) and applying a diagnonal method for $\delta
\downarrow 0$ and $T \uparrow \infty ,$ we see that, for any sequence $%
\varepsilon _n \to 0 ,$ there is a subsequence $V_{\varepsilon_n} $ (denoted
still by itself) such that $V_{\varepsilon_n}\longrightarrow V$ weakly in $%
H^1_{loc}(\Omega (\omega ))$. Moreover, (\ref{eq3.21}), (\ref{eq3.2}) and
Lebesque's domainated convergence theorem imply that
\begin{equation}  \label{eq3.22}
|V(x,t)|=1 \ \ a.e. \ \ in \ \ \Omega (\omega ).
\end{equation}
By taking the wedge product of $V_{\varepsilon _n} $ with the equation in (%
\ref{eq1.1}), we have
\[
V_{\varepsilon _n}\wedge \frac{\partial V_{\varepsilon _n}}{\partial t}=
div(V_{\varepsilon _n}\wedge \nabla V_{\varepsilon _n})+V_{\varepsilon
_n}\wedge(\nabla \omega \nabla V_{\varepsilon _n} ).
\]
Passing to the limit, we conclude that
\[
V \wedge \frac{\partial V}{\partial t}= div( V\wedge \nabla
V)+V\wedge(\nabla \omega \nabla V) \ \ in \ \ D^\prime (\Omega (\omega )).
\]
But (\ref{eq3.22}) yields
\[
V\frac{\partial V}{\partial t}= 0, \ \ \ div(V \nabla V)=0 \ \ in \ \
D^\prime (\Omega (\omega )) .
\]
Combing the last three equations with (\ref{eq3.22}), one gets

\[
\left \{
\begin{array}{cc}
V_t=\Delta V+V|\nabla V|^2 +\nabla \omega \nabla V & \ \ in \ \
D^{\prime}(\Omega (\omega )) \\
|V|=1 & \ \ in \ \ \Omega (\omega ) \\
u=g & \ \ on \ \ \partial \Omega .
\end{array}
\right.
\]
This completes the proof of Theorem 1.

 \section{Proof of Theorem 2}

 \setcounter{equation}{00}

 We assume $N=2$ as well as $(A_1), (A_2), (A_3)$ and $(A_4)$
throughout this section. Let $V_\var $ be classical solution to problem (\ref
{eq1.1}) . Then all conlusions in the last section hold true. We will
use covering arguments and ellitic estimates to prove
Theorem 2 by the coming  lemmas.

\noindent {\bf Lemma 10}\ \
  {\it For any $T_0> 0 ,$ there exist constants $C(T_0)>0$
and $\sigma _1 =\sigma _1(T_0)$, $\sigma _1 \in (0, \sigma (T_0)/4)$ with
the same $\sigma (T_0)$ as in Lemma 9 such that for all $t\in [0,
T_0], $ \ \ for all $x_0 \in \Omega \backslash\cup_{j=1}^m B_{\sigma
_1}(y_j(t))$ , for all $r\in (0, \sigma _1)$ and for all $\varepsilon \in
(0, \varepsilon _0)$, one has
\begin{eqnarray}  \label{eq4.1}
&\ \ & \int_{\partial \Omega_r (x_0)}|x-x_0||\frac{\partial V_\var}{\partial
\nu}|^2ds +\int_{\Omega _r(x_0)} \frac{1}{\varepsilon^2}(1-|V_\var|^2)^2dx
\nonumber \\
& \leq C(T_0) & \{ \int_{\partial \Omega _r(x_0)}|x-x_0|[|\frac{\partial
V_\var}{\partial T}|^2 +\frac{1}{4\varepsilon^2}(1-|V_\var|^2)^2]ds
\nonumber \\
& \ \ + &\int_{\Omega _r(x_0)}|x-x_0|[1+|\nabla V_\var |^2+|\nabla V_\var
||\partial _t V_\var |]dx\},
\end{eqnarray}
where $\Omega _r(x_0)=B_r(x_0)\cap \Omega $, $\nu$ and $T$ are,
respectively, the exterior unit normal vector and tangent vector of $
\partial \Omega _r$ such that $(\nu, T)$ is direct. }

{\bf Proof.}  By Lemma 9 and the fact that $\Omega $
is smooth, we can find constant $\alpha = \alpha (T_0), $ $\sigma _1=\sigma
_1 (T_0)\in (0, \sigma (T_0)/4)$ such that for all $t\in [0, T_0], $ all $%
x_0 =x_0(t)\in \Omega \backslash \cup_{j=1}^m B_{\sigma _1}(y_j (t))$ and
all $r\in (0, \sigma _1),$
\begin{equation}  \label{eq4.2}
(x-x_0)\cdot \nu\geq\alpha |x-x_0|, \forall \in \partial \Omega _r(x_0).
\end{equation}
Multiplying the equation in (\ref{eq1.1}) by $\nabla V_\var\cdot(x-x_0)$ and
integrate it over $\Omega _r = \Omega _r(x_0) .$ Neglecting the subscript $%
\varepsilon $, we obtain that
\begin{eqnarray}  \label{eq4.3}
& & \frac{1}{4\varepsilon ^2}\int _{\Omega _r} (1-|V|^2)^2 div (B(x)\cdot
(x-x_0))  \nonumber \\
&=& \frac{1}{4\varepsilon ^2}\int _{\partial \Omega _r} (1-|V|^2)^2 B(x)
\nu\cdot (x-x_0)  \nonumber \\
& & +I_3+I_4+I_5,
\end{eqnarray}
where
\[
I_3 =\displaystyle{\int_{\partial\Omega _r}}\displaystyle{\frac{\partial ( -
V)}{\partial \nu}}\cdot(\nabla V \cdot (x-x_0)) ,
\]
\[
I_4= \displaystyle{\int_{\Omega _r}}\nabla V\cdot\nabla(\nabla V(x-x_0)),
\]
and
\begin{eqnarray*}
I_5 &=& \displaystyle{\int_{\Omega _r}}\partial _t V\cdot(\nabla V\cdot
(x-x_0)) -(\nabla \omega \nabla V +A V)(\nabla V \cdot (x-x_0)) \\
& \leq & \int _{\Omega _r}C (|V|^2+|\nabla V|^2+|\partial _t V||\nabla
V|)|x-x_0|dx.
\end{eqnarray*}
By virtue of $(A_2)$ and the smallness of $r$ we may assume
\begin{equation}  \label{eq4.4}
div(B\cdot (x-x_0))\geq \lambda >0
\end{equation}
for all $x\in \Omega _r$ and some constant $\lambda $ depending only on $A$.
On the other hand, the integrand in $I_3$ can be writted as
\begin{equation}  \label{eq4.5}
- \frac{\partial V}{\partial \nu}\cdot[(\frac{\partial V}{\partial \nu}\nu+
\frac{\partial V}{\partial T}T)(x-x_0)] \leq \frac{-4}{5}|\frac{\partial V}{%
\partial \nu}|^2\nu\cdot (x-x_0) +C|\frac{\partial V}{\partial T}|^2|x-x_0|.
\end{equation}
The integrand in $I_4$ is nothing but $div (\frac{1}{2}|\nabla V|^2(x-x_0)).$
Hence, we have
\begin{equation}  \label{eq4.6}
I_4 = \int_{\partial \Omega _r}(\frac{1}{2}|\nabla V|^2\nu\cdot (x-x_0))ds
\end{equation}
Combing (\ref{eq4.2})-(\ref{eq4.6}) and using (\ref{eq3.1}), we have deduced
the desired (\ref{eq4.1}).

 \vskip 0.3cm

\noindent{\bf Lemma 11}\ \
 {\it For any interval $I\subset (0, \infty)$ with $|I|>0$ and any $\delta
 \in (0, \frac{1}{4}),$
there exist  $t\in I$ and $\eta_0>0$ such that
\begin{equation}  \label{eq4.7}
\int_{\Omega _{\delta /4}(t)}|\frac{\partial V_\var}{\partial t}|^2dx +
\int_{\Omega _{\delta /4} (t)}[ |\nabla V_\var|^2 + \frac{1}{\varepsilon^2}%
(1-|V_\var|^2)^2]dx \leq C(t,\delta )< +\infty
 \end{equation}
for all $\varepsilon \in (0,
\eta _0) $. Moreover,
there exist a $\varepsilon _1 \in (0, \eta _0)$ such
that
\begin{equation}  \label{eq4.8}
\{x\in \bar \Omega : |V_\var(x,t)|< \frac{1}{2}\}\subset \cup_{j=1}^m
B_{\delta /4}(y_j(t)), \forall \varepsilon \in (0, \varepsilon_1).
\end{equation}
}

{\bf Proof.} The proof is almost the same as the one of Lemma 4.2 in [18].
Here we copy it just for the convenience.

First, Using (3.21) and standard methods in real analysis, we easily get
(4.7). Thus, what we need do is only to prove (\ref{eq4.8})
for sufficiently small $\delta >0$ and the same $t$ as in (\ref{eq4.7}). If
the conclusion were not true, we could find $\delta_1 \in(0, \sigma (t))$
(with the same $\sigma (t)$ as in Lemma 9), a sequence $%
\varepsilon_k\searrow 0$, $\varepsilon_k\in(0, \delta_1)$, and $\{x_k\}
\subset\bar{\Omega}\setminus \cup_{j=1}^mB_{\delta_1}(y_j(t))$ such that $%
|V_{\varepsilon_k}(x_k ,t)| <\frac{1}{2}$ for all $k$. Hence, by virtue of (%
\ref{eq3.2}) and the fact $|V_\var|= |g|= 1$ on $\partial \Omega$, we see
that there is a ball $B_{C_1\varepsilon_k}(x_k)\subset \Omega\setminus
\cup_{j=1}^m B_{\frac{\delta_1}{2}} (y_j(t))$ for some constant $C_1>0$ with
$|V_{\varepsilon_k}(x, t)|\leq \frac{3}{4}$ for all $x\in
B_{C_1\varepsilon_k}(x_k)$ and all sufficiently large $k\geq k_0.$ Let $r_k
=C_1\varepsilon_k$, $B_k =B_{r_k}(x_k)$ and $V_k=V_{\varepsilon_k}$. Since $%
N=2 ,$ one has that
\begin{equation}  \label{eq4.9}
\varepsilon_k^{-2}\int_{B_k}(1-|V_k(x,t)|^2)^2dx \geq C_2>0
\end{equation}
for all $k\geq k_0$ and some positive constant $C_2$ depending only on $C_1$%
.

 On the other hand, as $r_k \to 0, $ (\ref{eq4.7}) implies that
\begin{equation}  \label{eq4.10}
\int_{\Omega _{\sqrt{r_k}}(x_k)}|\frac {\partial V_k}{\partial t}|^2dx+
\int_{\Omega _{\sqrt{r_k}}(x_k)} [|\nabla V_k|^2 + \varepsilon_k^{-2}
(1-|V_k|^2)^2 ]dx\leq C(t, \delta )
\end{equation}
for all $k\geq k_0$. Thus, letting
\[
f_k(r)=\int_{\partial B_r(x_k)\cap \Omega }[ |\nabla V_k|^2 +
\varepsilon_k^{-2}(1-|V_k|^2)^2 ]ds
\]
we have
\begin{eqnarray*}
C(t, \delta) &\geq & \displaystyle{\int}_{(B_{\sqrt{r_k}}(x_k)\setminus
B_{r_k}(x_k)) \cap \Omega}[\nabla V_k|^2 + \varepsilon_k^{-2}
(1-|V_k|^2)^2]dx \\
&=& \displaystyle{\int}_{r_k}^{\sqrt{r_k}}f_k(r)dr \\
&\geq& \frac{1}{2}|\ln r_k|\min_{r_k\leq r\leq\sqrt{r_k}} \{rf_k(r)\}, \ \
\forall k\geq k_0.
\end{eqnarray*}
Therefore, for each $k\geq k_0,$ we can find a $\lambda _k \in(r_k, \sqrt{r_k%
})$ such that
\begin{equation}  \label{eq4.11}
\lambda_k f_k(\lambda_k) \leq 2|\ln \lambda _k|^{-1}C(t, \delta ).
\end{equation}
Using Lemma 10 for $x_0=x_k$ and $r=\lambda_k$, we obtain
\begin{eqnarray*}
\varepsilon_k^{-2}\displaystyle{\int}_{B_k}(1-|V_k|^2)^2 &\leq &
\varepsilon_k^{-2}\displaystyle{\int}_{\Omega_{\lambda_k} (x_k)}(1-|V_k|^2)^2
\\
&\leq & C\lambda _k [ \displaystyle{\int}_{\Omega _{\lambda
_k}(x_k)}(1+|\nabla V_k|^2+|\nabla V_k| |\partial _t V_k|)dx \\
&+& \displaystyle{\int}_{\partial \Omega }|\frac {\partial V_k}{\partial
T}|^2ds + \displaystyle{\int}_{\partial B_{\lambda _k}\cap\Omega}(|\frac
{\partial V_k}{\partial T}|^2 +\displaystyle{\frac{(1-|V_k|^2)^2}{%
4\varepsilon^2_k}})ds] \\
&\leq & C [ \lambda _k +2|\ln \lambda _k|^{-1}] \ \ \ (by \ \ (\ref{eq4.10}%
)\ \ and \ \ (\ref{eq4.11})).
\end{eqnarray*}
This contradicts with (\ref{eq4.9}) because of the fact $\lambda_k \to 0 .$
In this way, we finish the proof of Lemma 11.

{\it {\rm \vskip 0.3cm }}

\noindent{\bf Lemma 12 }\ \
  {\it With the same $\delta , t $ and $\varepsilon _1$ as
in Lemma 11, one has that for all $x\in\bar{\Omega}\setminus
\cup_{j=1}^mB_\dd (y_j(t)),$
\begin{equation}  \label{4.12}
1-C(\delta , t)\varepsilon^{\frac{1}{2}} \leq |V_\var (x, t)|\leq 1+C(\delta ,
t)\varepsilon^{\frac{1}{2}} , \ \ \forall \varepsilon\in(0, \varepsilon_1).
\end{equation}     }

{\bf Proof. } we follow  the proof of Lemma 4.3 in [18]. Here we give the details just for the reader's convenience.
Fix $x_0\in \partial \Omega $ and let $\sigma _2
=\frac{1}{4} \min \{ 1 , [ \sigma (T)]^2 \}.$ Using (\ref{eq4.7}) and the
arguements from (\ref{eq4.10}) to (\ref{eq4.11}) one can easily see that
\begin{equation}  \label{eq4.13}
\int_{\partial B_{\lambda _\var}(x_0)\cap \Omega }[|\nabla V_\var (x,t)|^2+
\frac{(1-|V_\var (x,t)|^2)^2} {\varepsilon ^2} ]dx\leq C(\delta , t, \sigma
_2)
\end{equation}
for some $\lambda _\var \in [\sigma _2 , \sqrt{\sigma _2}]$ and all $%
\varepsilon \in (0, \varepsilon _1)$. Moreover, (\ref{eq4.13}), (\ref{eq4.7}%
) and Lemma 10 imply that $\int_{\partial \Omega }|\frac {\partial
V_\var }{\partial N}|^2 \leq C $ independent of $\varepsilon $. Therefore,
we deduce that
\begin{equation}  \label{eq4.14}
\int_{\partial\Omega}|\nabla V_\var|^2 ds+ \int_{\Omega _{\delta /4}(t)}
(|\nabla V_\var|^2+ |\frac{\partial V_\var }{\partial t}|^2+\frac{%
(1-|V_\var|^2)^2 }{\varepsilon^2})dx \leq C(\delta, t)
\end{equation}
and
\begin{equation}  \label{eq4.15}
\alpha _1 \geq |V_\var(x,t)|^2\ge \alpha _2>0 , \,\,\,\,\,\,\,\,\forall x\in
\Omega _{\delta /4}(t) \ \ (by \,\,\,\, (\ref{eq3.2}) \,\,\,\, and \,\,\,\, (%
\ref{eq4.8}))
\end{equation}
for all $\varepsilon\in(0, \varepsilon_1)$. Let $R_0$ be a positive constant
to be determined later. As we will see, it depends only on $\alpha_1$, $%
\alpha_2$ and $C(\delta ,t )$ in (\ref{eq4.14}). Fix a constant $r_0\in (0,
\min\{\frac{\delta}{4}, R_0\})$ which will be suitably small at last. For an
arbitrary $y\in\Omega _\dd(t)$, choose a number $R\in (0, \min\{\frac{\delta%
}{8}, \frac{R_0-r_0}{2}\})$ satisfying $B_{2R+r_0}(y)\subset\Omega _{\frac{%
\delta}{2}}(t)$.

 Write
\[
V_\var(x , t) =\rho_\var(x, t)e^{i\psi_\var(x,t)}
\]
on $B_{2R+r_0}(y)\times (t-C_\var , t+C_\var)$ for some $C_\var >0$(see (\ref
{eq3.2})) so that the equation in (\ref{eq1.1}) turns to be
\begin{equation}  \label{eq4.16}
div(\rho_\var^2\nabla\psi_\var) =\rho_\var^2(\partial _t \psi_\var
-\nabla \psi_\var \nabla \omega )\,\,\,\,\,\,\,\,
\mbox{in
$B_{2R+r_0}(y)$}
\end{equation}
and
\begin{equation}  \label{eq4.17}
\Delta\rho_\var + B\frac{(1-\rho_\var^2)}{\varepsilon^2}\rho_\var =
|\nabla\psi_\var|^2\rho_\var+\partial _t \rho_\var -\nabla \rho _\var \nabla
\omega -A \rho _\var \,\,\,\,\,\,\,\,\mbox{in $B_{2R+r_0}(y)$}.
\end{equation}
Moreover, using (\ref{eq4.14}) and the Fubini's theorem(see the arguements
from (\ref{eq4.10})to the (\ref{eq4.11})), one can find $R_\var\in(R, R+%
\frac{r_0}{2})$ such that
\begin{equation}  \label{eq4.18}
\int_{\partial B_{R_\var}(y)} |\nabla V_\var|^2+\frac{(1-|V_\var|^2)^2 }{%
\varepsilon^2} \leq 2r_0^{-1}C(\delta , t).
\end{equation}
It easily follows from (\ref{eq4.18}) and Lemma 7 that
\[
\max_{x\in\partial B_{R_\var}(y)}|1-|V_\var(x,t)|^2|\leq C(r_0^{-1}, \delta
)\varepsilon^{\frac{1}{2}}
\]
which, together with (\ref{eq3.1}), implies that
\begin{equation}  \label{eq4.19}
|1-|V_\var(x,t)||\leq C_2\varepsilon^{\frac{1}{2}}
\end{equation}
for all $x\in\partial B_{R_\var}(y)$ , all $\varepsilon\in(0, \varepsilon_1)$
and some constant $C_2$ depending only on $r_0^{-1}$, $\delta$ and $t$.

On one hand, applying Theorem 2.2 of Chapter V in \cite{g} to the
equation (\ref{eq4.16}) for $\psi_\var$ with the coefficient $\rho_\var$
satisfying (\ref{eq4.15})and using the notation $\oint_E=\frac{1}{E}\int_E$,
we obtain that
\begin{eqnarray}  \label{eq4.20}
(\oint_{B_{R+\frac{r_0}{2}}(y)}|\nabla\psi_\var|^pdx)^{\frac{1}{p}} & \leq &
C_3\{(\oint_{B_{2R+r_0}(y)}|\nabla\psi_\var|^2dx)^{\frac{1}{2}}  \nonumber \\
& + & R[\oint_{B_{2R+r_0}(y)}(|\nabla \psi||\nabla \omega |+|\partial
_t\psi_\var |) ^{\frac{2p}{3}}dx]^{\frac{3}{2p}}\}
\end{eqnarray}
for some $p\in (2, 3]$ depending only on $\alpha_1$, $\alpha _2$ in (\ref
{eq4.15}), for some $R_0 >0$ depending only on $\alpha _1$, $\alpha _2$ and $%
C$ in (\ref{eq4.14}), for some $C_3>0$ depending only on $p$, $\alpha _1$
and $\alpha _2$, for all $R<\frac{R_0-r_0}{2}$, and for all $%
\varepsilon\in(0, \varepsilon_1)$.

On the other hand, equation (\ref{eq4.17}) implies that the
function
\[
\bar{U}_\var(x,t)\equiv 1 - \rho_\var(x,t)
\]
satisfies
\[
-\Delta \bar{U}_\var + \frac{C_\var(x)}{\varepsilon^2}\bar{U}_\var = f_\var
\]
in $B_{R_\var}(y)$ with $0< C(\alpha_1)\leq C_\var(x)\equiv
B(1+\rho_\var)\rho_\var\leq C(\alpha_2)$ and
\[
f_\var\equiv |\nabla\psi_\var|^2\rho_\var +\partial _t\rho_\var -\nabla
\omega \nabla \rho_\var-A\rho_\var \in L^{\frac{p}{2}} (B_{R_\var}(y)) .
\]
By (\ref{eq4.20}) and (\ref{eq4.14}), we see that
\[
||f_\var||_{L^{\frac{p}{2}} (B_{R_\var}(y))}\leq C(\alpha_1, \alpha_2,
\delta) ,\ \ \ \ \forall \varepsilon\in(0, \varepsilon_1).
\]
Moreover, (\ref{eq4.19}) yields $-C_4\varepsilon^{\frac{1}{2}}\leq\bar{U}%
_\var \leq C_4 \varepsilon^{\frac{1}{2}}$ on $\partial B_{R_\var}(y). $
Therefore, a standard elliptic estimate (Theorem 8.16 in \cite{gt}) gives us
\[
|\bar{U}_\var|\leq C(C_2, \alpha_1, \alpha_2)\varepsilon^{\frac{1}{2}} \ \
in \ \ B_{R_\var}(y)\times \{ t \}.
\]
Particularly, we have
\begin{equation}  \label{eq4.21}
1-C\varepsilon^{\frac{1}{2}}\leq \rho_\var(x) = |V_\var(x,t)|\leq 1+
C\varepsilon^{\frac{1}{2}} ,\ \ \forall x\in B_R(y), \ \ \forall \varepsilon
\in (0, \varepsilon _1).
\end{equation}

 Now for any $G_0\subset\subset\Omega$, choose $r_0 = \frac{1}{4}%
\min\{ R_0, \delta, dist(G_0, \partial \Omega)\}$ and fix
\[
R=\min\{\frac{\delta}{8}, \frac{R_0-r_0}{2}, \frac{1}{4}dist(G_0,
\partial\Omega)\}.
\]
Then, by the arbitrariness of $y\in \bar{G_0}\setminus \cup_{j=1}^m B_\dd
(y_j(t))$, we can find finite balls, $B_R(y_i)$, $i=1, 2, \cdots, n$, such
that $\cup^n_{i=1}B_R(y_i)\supset \bar{G}_0\setminus
\cup_{j=1}^mB_\dd(y_j(t))$ and (\ref{eq4.21}) holds true for all $x\in
B_R(y_i)$ and all $i=1, 2, \cdots, n$. In this way, we conclude
\begin{equation}  \label{eq4.22}
1-C\varepsilon^{\frac{1}{2}} \leq |V_\var(x,t)|\leq 1+C\varepsilon^{\frac{1}{%
2}},\,\,\,\,\forall x\in \bar{G_0}\setminus \cup_{j=1}^mB_\dd(y_j(t)),
\forall \varepsilon \in (0, \varepsilon _1 )
\end{equation}
for some constant $C$ and $\varepsilon_1$ both independent of $\varepsilon$.
Moreover, using the fact $\int_{\partial\Omega}|\nabla V_\var|^2\leq
C(\delta)$(see (\ref{eq4.14})) and repeating the argument above for $%
\Omega_R = B_R(y)\cap\Omega$ with $y\in\partial\Omega$, we can find a domain
$G^\prime\subset\subset\Omega$ such that (\ref{eq4.22}) holds true for all $%
x\in\bar{\Omega}\setminus G^\prime$. Combining this result with (\ref{eq4.22}%
), we have proved Lemma 12.

Now, combining Lemmas 11 and 12,
 we have completed the proofs of conclusion (i) of Theorem 2
 as well as the corresponding part of Corollary 3.
\vskip 0.3cm

Next, we are going to prove the second part of Theorem 2.

\noindent {\bf Lemma 13}\ \
{\it $\{\psi _\varepsilon : \varepsilon\in(0, \varepsilon _1)\}$
is compact in $H^1(B)$ with $B=B_{R+\frac{r_0}{2}}(y).$}

{\bf Proof.} (4.14) and (4.15) imply that
$\{ \nabla \rho_\var ^2   : \varepsilon\in(0, \varepsilon _1)\}$
is bounded in $L^2(B) .$ Thus, by (4.20) and Holder inequality, we
see that
$\{ \nabla \rho_\var ^2\nabla \psi _\var
: \varepsilon\in(0, \varepsilon _1)\}$  is bounded in $L^q(B)$
with $q=\frac{2p}{p+1}>1.$ Hence
$\{  \psi _\var
: \varepsilon\in(0, \varepsilon _1)\}$  is bounded in $W^{2,q}(B)$
by (4.7), (4.15), equation (4.16) and standard elliptic estimates.
Furthermore, the Rellich-Kondrachov's theorem tells us that the set
$\{  \psi _\var
: \varepsilon\in(0, \varepsilon _1)\}$  is compact in $H^1(B) .$

\vskip 0.3cm

To complete the proof of conclusion (ii) of Theorem 2, we need
only to prove estimate (1.6). This is because the strong convergence of
$V_{\varepsilon_n} \to V$ is a direct consequence of Theorem 1,
Lemma 13, and (1.6).

In order to prove (1.6), we rewrite (4.17) as
\[\Delta\rho_\var + \frac{(1+\var ^2 -\rho_\var^2)}{\varepsilon^2}\rho_\var =
|\nabla\psi_\var|^2\rho_\var+\partial _t \rho_\var -\nabla \rho _\var \nabla
\omega -(A-1) \rho _\var .\]
Multiplying this equality by $\rho_\var $ and integrating it over
$B_{R_\var}(y), $ we have

\begin{eqnarray*}
\int_{B_{R_\var}(y)} \frac{(1+\var ^2 -\rho_\var^2)}{\varepsilon^2}
\rho_\var ^2 dx & = &\int_{B_{R_\var}(y)}  [ |\nabla\rho_\var |^2+
|\nabla\psi_\var|^2\rho_\var ^2+(\partial _t \rho_\var -\nabla
\rho _\var \nabla \omega \\
& - & (A-1) \rho _\var )\rho_\var ]dx
- \int_{\partial B_{R_\var}(y)}\frac{\partial \rho_\var}{\partial \nu}
\rho_\var .
\end{eqnarray*}

It is easy to see that the sum on the right side hand can be bounded by a
constant $C_3=C(\alpha_1, \alpha_2, \delta, t, r_0^{-1})$ due to (4.7),
(4.15) and (4.18). Theorefore, applying (4.15), corollary 8 and conclusion (i)
of Theorem 2, we obtain that
\begin{eqnarray*}
\int_{B_{R_\var}(y)} \frac{(1+\var ^2 -\rho_\var^2)^2}{\varepsilon^2}dx
& \leq & \alpha_2 ^{-2}\int_{B_{R_\var}(y)}
\frac{(1+\var ^2 -\rho_\var^2)^2}{\varepsilon^2}\rho_\var ^2 dx \\
&\leq & \max_{(x,t)\in \Omega (\omega_t^{\delta/2})}(1+\var ^2 -\rho_\var ^2)
C_3 \\
&\leq & C(C_3, \delta , t)\var^{\frac{1}{2}}.
\end{eqnarray*}
Hence
\begin{eqnarray}               \label{eq4.23}
\var ^{-2} \int_{B_{R_\var}(y)} (1-|V_\var |^2)^2dx
& = & \var ^{-2}\int_{B_{R_\var}(y)}
[(1+\var ^2 -\rho_\var^2)^2-\var ^4 -2\var ^2 (1-\rho_\var ^2)]dx \nonumber \\
 &\leq & C \var^{\frac{1}{2}}.
\end{eqnarray}

On the other hand, setting $P_{\var }=1-|V_\var |^2,$ we have, by equation
(4.17), that
$$
-\Delta P_\var + 2\frac{\rho_\var^2 P_\var}{\varepsilon^2} =
2|\nabla\psi_\var|^2\rho_\var ^2+\partial _t \rho_\var ^2
+2|\nabla \rho_\var |^2 -\nabla \rho _\var ^2 \nabla
\omega -2A \rho _\var ^2$$
in $B_{R_\var }(y)$.
Multiplying this equation by $P_\var  ,$ integrating the resulted
equality over   $B_{R_\var }(y),$ and using (4.7), (4.15), (4.18),
(4.23), the boundedness of
$|\psi_\var |^2$ in $L^2(B) ,$  and the estimate
$$|P_\var |=(1+|V_\var |)(1-|V_\var |)\leq C(\delta , t)\var ^{\frac{1}{2}}$$
from the conclusion (i) of Theorem 2, we obtain that
\begin{equation} \label{eq4.24}
\int_{ B_{R_\var }(y)}|\nabla P_\var |^2 dx\leq
 C(\delta , t)\var ^{\frac{1}{2}}.
 \end{equation}
 Observing  $\nabla |V_\var |=\nabla \rho_\var = -(2\rho_\var )^{-1}
 \nabla P_\var ,$ we have, by (4.15) and (4.24), that
\begin{equation} \label{eq4.25}
\int_{ B_{R_\var }(y)}|\nabla |V_\var | |^2 dx\leq
 C(\delta , t)\var ^{\frac{1}{2}}.
 \end{equation}
 Finally, using (4.23), (4.25) and Lemma 13, and repeating the covering
 arguments in the end of the proof of (i) of Theorem 2, we have proved the
 (1.6) and thus
 (ii) of Theorem 2.


\begin{thebibliography}{99}
\bibitem{cr}   S. J. Chapman and G. Richardson, Vortex pinning by
inhomogeneities in type-II superconductors, {\em Phys. D.}, vol. 108(1997),
397-407.

\bibitem{cdg}   S. J. Chapman, Q. Du and M. D. Gunzburger, A
Ginzburg-Landau type model of superconducting/ normal junctions including
Josephson junctions, {\em Euro. J. Appl. Math.}, vol. 6(1995), 97-114.

\bibitem{r}  J. Rubinstein, On the equilibrium position of
Ginzburg Landau vortices, {\em Z. Angew. Math. Phy.}, 46(1995), 739-751.

\bibitem{as}   N. Andr\'{e} and I. Shafrir, Asymptotic behaviour
for the Ginzburg-Landau functional with weight I, II, {\em Arch. Rat. Mech.
Anal.}, 142(1998), 45-98.

\bibitem{dl}   S. Ding and Z. Liu, Remarks for the asymptotics of a
kind of Ginzburg-Landau functional (in chinese), {\em Chin. Ann. Math.},
19A(1998), 621-628.

\bibitem{j1}   H.Y. Jian, A relation between $\Gamma$-convergence
of functionals and their associated gradient flows, {\em Science in China
(Series A)}, 42(2)(1999), 133-139.

\bibitem{jw1}   H.Y. Jian and Y.D. Wang, Ginzburg-Landau vortices
in inhomogeneous superconductor, {\em
Advance Mathematic}, 28(1)( 1999), 83-84.

\bibitem{dg}   Q. Du and M.D. Gunzburger, A model for
superconducting thin films having variable thickness, {\em Physica D.},
69(1993), 215-231.

\bibitem{cdg1}  S. J. Chapman, Q. Du and Gunzburger, A model for
variable thickness superconducting thin films, {\em Z. Angew Math. Phys.},
47(1996), 410-431.

\bibitem{jw2}   H.Y. Jian and Y.D. Wang, Ginzburg- Landau vortices
with pinning functions and self-similar solutions in harmonic maps,
{\em Science in China Ser. A, } 43(2000), to appear.

\bibitem{bbh}   F. Bethuel, H. Brezis and F. H\'{e}lein,
Ginzburg-Landau vortices, Birkh\"{a}user, Boston , 1994.

\bibitem{s}   M. Struwe, On the asymptotic behaviour of minimizers
of the Ginzburg-Landau model in 2-dimensions, {\em J. Diff. Int. Equations},
7(1994), 1613-1624.

\bibitem{l2}   F. H. Lin, Some dynamical properties of
Ginzburg-Landau vortices, I, II, {\em Commu. Pure Appl. Math.}, vol.
49(1996), 323-364.


\bibitem{js1}   R.L. Jerrard and H.M. Soner, Scaling limits and regularity
for a class of  Ginzburg-Landau systems, {\em Ann. Inst. Henri. Poincare
Analyse
Nonlineaire. }, 16(1998), 423-446.



\bibitem{js}   R.L. Jerrard and H.M. Soner, Dynamics of Ginzburg
Landau vortices, {\em Arch. Rat. Mech.}, 142(1998), 99-125.

\bibitem{l3}   F.H. Lin, Complex Ginzburg-Landau equations and
dynamics of vortices, filaments, and codimension-2 submanifolds, {\em Commu.
Pure Appl. Math.}, vol. 51(1998), 385-441.

\bibitem{j2}   H.Y. Jian, The dynamical law of Ginzburg-Landau
vortices with a pinning effect, {\em Appl. Math. Letters}, 13(2000), 91-94.

\bibitem{j3}    H.Y. Jian  and B. Song, Vortex dynamics of
Ginzburg-Landau equations in inhomogeneous superconductors,
 {\em  Journal of  Differential Equations}, 170(2001),123-141.

 \bibitem{zd}  Z. F. Zhang and T. R. Ding, Stability theorey of
differential equations (in Chinese), Science Press, Beijing, 1997.

 \bibitem{g}  M. Giaquinta, Multiple integrals in the calculus of
variations and nonlinear elliptic systems, Princeton Univ. Press, Princeton,
1983.

\bibitem{gt}  D. Gilbarg and N.S. Trudinger, Elliptic partial
differential equations of second order, Springer-Verlag, Heidelberg, New
York, 1977.
\end{thebibliography}
\end{document}